\documentclass[12pt]{article}

\usepackage[utf8]{inputenc}
\usepackage[margin=3cm]{geometry}
\usepackage{amsfonts, amsmath, amssymb, amsthm}
\usepackage{indentfirst}
\usepackage{xcolor}

\newtheorem{Dfn}{Definition}
\newtheorem{Thm}{Theorem}
\newtheorem{Prp}{Proposition}
\newtheorem{Lem}{Lemma}

\newenvironment{Prf}{\par\noindent{\it Proof:}\rm}{\newline\rightline{\textbf{QED}}\par}
\newtheorem{Cor}{Corrolary}

\newtheorem{Prb}{Problem}

\newcommand{\R}{\mathbb{R}}

\newcommand{\con}{\mathfrak{c}}

\newcommand{\Cl}{\rm{Cl}}

\title{
    \textbf{Application of fusion technique to the solution of Harrington problem} \textit{and its generalizations to Baire functions}, \textbf{part I}
}

\author{Sławomir Kusiński}

\date{}

\begin{document}
    \maketitle
    
    \begin{abstract}
        In this paper we provide solutions of the Harrington problem (along with a few generalizations) proposed in a book \textit{Analytic Sets}. The original problem asks if for arbitrary sequence of continuous functions from \( \R^\omega \) to a fixed compact interval we can find a subsequence point-wise convergent on some product of perfect subsets of \( \R \). We reduce aforementioned problem to functions from \( C^\omega \) to \(C\), where \(C\) is a standard Cantor set as well as also provide solution to the problem with Baire functions in place of continuous ones. Our main focus is on showing applications of the fusion lemma - a result about perfect trees used among others to prove minimality of Sack's forcing - to the problem at hand.
    \end{abstract}
    
    \section{Preliminary notions}
    
    By a perfect set we mean a subset of a topological space that is compact and has no isolated points. It is worth noting that some authors use a slightly weaker definition of a perfect set, namely they require it to be only closed instead of compact. By \( 2 \) we mean a set \( \{ 0, 1 \} \) with discrete topology. By \( C \subseteq \R \) we mean a standard Cantor set. It is well known that \( C \) is homeomorphic to the space \( 2^\omega \) with product topology where \( \omega \) denotes the set of all natural numbers.
    
    In \cite{AS} in the problems section L Harrington published a following problem (as a possible weakening of a problem of Halpern).
    \begin{Prb}
        Given continuous functions \( f_n \colon \R^\omega \to [0; 1] \), do there exist a set \( N = \{ n_i \colon i \in \omega \} \in [\omega]^\omega \) and nonempty perfect sets \( P_j \subseteq [0;1] \) for \( j \in \omega \) such that the subsequence \( (f_{n_i})_{i \in \omega} \) is pointwise convergent on the product \( \prod\limits_{j \in \omega}  P_j \)?
    \end{Prb}
    
    We have acquired information that the one dimensional version of such problem has been solved in 1920s by S Mazurkiewicz, but unfortunately we were not able to trace it back to the original paper. In \cite{Lav} Laver showed amongst others that we get an equivalent problem if we substitute continuous for measurable functions or functions with the Baire property. This prompted to consider variants of the problem with different notions of measurability and different topologies (most notably with Ellentuck topology \cite{El}).
    
    We will restrict our attention to the functions with domain \( C^\omega \) and codomain \( C \). Note that restriction of domain is not a weakening of the statement and as we will show such a restriction of codomain leads in fact to an equivalent problem, ie by solving such a variant we will also solve the original problem in all its generality. Apart from the classical variant for continuous function we will also consider Baire functions and measurable functions. Our main tool at work will be the fusion lemma - a result about perfect trees, that originated in forcing theory. In its original form it was used to prove the minimality of Sack's forcing. We will use a slightly more complex, topological variant of it.
    
    \begin{Dfn}
        Let \( T \subseteq 2^{<\omega} = \bigcup\limits_{n \in \omega} 2^n \).
        
        We will say that \( T \) is a tree if for any \( s \in T \) and \( t \in 2^{<\omega} \) if \( t \subseteq s \) then \( t \in T \), ie it is closed under taking the initial segment.
        
        We will further say that \( T \) is a perfect tree if for any \( s \in T \) there exist \( t_1, t_2 \in T \) with \( t_1 \not= t_2 \) such that \( s \subset t_1 \) and \( s \subset t_2 \).
    \end{Dfn}
    
    \begin{Dfn}
        Let \( T \) be a perfect tree and for each \( s \in T \) let there be a perfect set \( U_s \subseteq X \), where \( X \) is a metrizable space. Then \( (U_s)_{s \in T} \) is called a fusion sequence if
        \begin{itemize}
            \item \( U_{s_1} \supseteq U_{s_2} \) for \( s_1, s_2 \in T \) and \( s_1 \subseteq s_2 \),
            \item \( U_{s ^\smallfrown  0} \cap U_{s ^\smallfrown  1} = \emptyset \) for \( s \in T \).
        \end{itemize}
    \end{Dfn}
    
    We have a following fundamental property of fusion sequences.
    
    \begin{Thm}
        (Fusion Lemma) Let \( T \) be a perfect tree and let \( \tilde{T} = \{ f \in 2^\omega \colon \forall_{n \in \omega} f|_n \in T \} \). Let \( (U_s)_{s \in T} \) be a fusion sequence. If the diameter of \( U_s \) tends to \( 0 \) with increasing length of \(s\) then the set
        \[
        P = \bigcap\limits_{n \in \omega} \bigcup\limits_{s \in T} U_s =  \bigcup\limits_{f \in \tilde{T}} \bigcap\limits_{n \in \omega} U_{f|n}
        \]
        is a perfect set and it is homeomorphic to the Cantor set.
        \cite{Jech}\cite{JechM}
    \end{Thm}

    If we are considering exclusively subsets of the Cantor set (or for that matter a space homeomorphic to it) the restriction on diameter can be dropped up to extent.
    
    \begin{Thm}
        (Fusion Lemma) Let \( T \) be a perfect tree and let \( \tilde{T} = \{ f \in 2^\omega \colon \forall_{n \in \omega} f|_n \in T \} \). Let \( (U_s)_{s \in T} \) be a fusion sequence of subsets of \(C\). Then the set
        \[
            Q = \bigcap\limits_{n \in \omega} \bigcup\limits_{s \in T} U_s =  \bigcup\limits_{f \in \tilde{T}} \bigcap\limits_{n \in \omega} U_{f|n}
        \]
        contains a non-empty perfect subset. If we further assume that the sets \( U_s \) are all basic clopen subsets of \(C\) (with respect to the product topology on \(2^\omega\)) then \( Q \) itself is perfect.
    \end{Thm}
   
   \section{Cantor set, perfect sets and perfect trees}
   
    It worth noting that the Cantor set, perfect sets and perfect trees are very closely connected. In this section we will outline those properties of them that will be useful in proving our main result. For any non-empty set \( A \subseteq  2^\omega = C \) let
    \[
        T_A = \{ x|_n \colon x \in A, n \in \omega \}.
    \]
    \(T_A\) is a perfect tree if and only if the set \(A\) has no isolated points. Moreover the set
    \[
        \tilde{T_A} = \{ x \in 2^\omega \colon \forall_{n \in \omega} x|_n \in T_A \} = \Cl(A),
    \]
    ie \( A = \tilde{T_A} \) if and only if \( A \) is perfect. On the other hand if \(T\) is any perfect tree the fusion lemma automatically gives us that \( \tilde{T} \) is a perfect set. Of course \( T_{\tilde{T}} = T \).
   
    \begin{Cor}
        Let \( P \subseteq C \) be perfect and non-empty. Then \( P \) is homeomorphic to \( C \).
    \end{Cor}

    \begin{Prp}
        Let \(X\) be a metric space and \( A \subseteq X \) be a non-empty perfect set. Then there exist a subset of \( A \) homeomorphic to the Cantor set.
    \end{Prp}
    \begin{Prf}
        As \( A \) is perfect and non-empty it has infinitely many elements. Let \( x_0, x_1 \in A \) be distinct. There has to exist \( r_0 > 0 \) such that \( B(x_0, r_0) \cap B(x_1, r_0) = \emptyset \). Let
        \begin{align*}
            P_{0} = A \cap \bar{B}(x_0, \frac{r_0}{2}),\\
            P_{1} = A \cap \bar{B}(x_1, \frac{r_0}{2}).
        \end{align*}
        Observe that as \(A\) is perfect so are the sets defined above. Now inductively for every \( i_0, \ldots, i_n \in 2 \) there exist distinct \( x_{i_0, \ldots,i_n,0}, x_{i_0, \ldots,i_n,1} \in P_{i_0,\ldots,i_n} \) and there exists \( r_{n+1} > 0 \) (common for all the sequences of length \(n+1\)) such that \( B(x_{i_0, \ldots,i_n,0}, r_{n+1}) \cap B(x_{i_0, \ldots,i_n,1}, r_{n+1}) = \emptyset \). We can thus define
        \begin{align*}
            P_{i_0, \ldots,i_n,0} = A \cap \bar{B}(x_{i_0, \ldots,i_n,0}, \frac{r_{n+1}}{2}),\\
            P_{i_0, \ldots,i_n,1} = A \cap \bar{B}(x_{i_0, \ldots,i_n,1}, \frac{r_{n+1}}{2}).
        \end{align*}
        It is clear from the construction that \( \lim\limits_{n \to +\infty} r_n = 0 \). From the compactness it follow that for any \( y \in 2^\omega \) we have \( |\bigcap\limits_{n \in \omega} P_{y|_n}| = 1 \). Thus the set
        \[
            P = \bigcup\limits_{y \in 2^\omega} \bigcap\limits_{n \in \omega} P_{y|_n}
        \]
        is homeomorphic to the Cantor set.
    \end{Prf}
    
    \begin{Cor}
        Let \( f \colon \R^\alpha \to \R \), where \( \alpha < \omega_1 \), be continuous and let \( A_n \subseteq \R \) for \( n \in \alpha \) be perfect. If \( f(\prod\limits_{n \in \alpha} A_n) \) contains a perfect subset then it contains a subset homeomorphic to the Cantor Set.
    \end{Cor}

    Topological spaces that do not have a non-empty dense-in-itself subset are called scattered. It is a well know property that second countable scattered spaces are countable.
    
    \begin{Prp}
         Let \(X\) be a compact metric space and \( f \colon X \to \R \) be continuous. If there exists \( B \subseteq f(X) \) homeomorphic to the Cantor set then there exists \( A \subseteq X \) also homeomorphic to the Cantor set and such that \( f(A) \subseteq B \). 
    \end{Prp}
    \begin{Prf}
        As \( f^{-1}(B) \) is compact, if it didn't have any perfect subsets then it would be in fact scattered as closure operation preserves isolated points. On the other hand \( f^{-1}(B) \) cannot be scattered as it is uncountable. It means that there exist a non-empty perfect subset \( P \subseteq f^{-1}(B) \) and from the corollary above we get that \( P \) has to have a subset homeomorphic to the Cantor set.
    \end{Prf}

    The next two propositions will be vital in some of our later arguments. 

    \begin{Prp}
        Let \( A \subseteq C \) be a dense \(G_\delta\) set. There exists a non-empty perfect set \( P \subseteq A \).
    \end{Prp}
    \begin{Prf}
        Let \( A = \bigcap\limits_{i \in \omega} U_i  \), where the sets \( U_i \) are open and dense. Let \( W_{(0)}, W_{(1)} \subseteq U_0 \) be disjoint basic clopen sets. Now with \( W_s \) defined for \( s \in 2^{i} \) let \( W_{s^\smallfrown 0}, W_{s^\smallfrown 1} \subseteq W_s \cap U_{i+1} \) be disjoint basic clopen sets. From the properties of fusion we obtain that the set
        \[
            P = \bigcup\limits_{x \in 2^\omega} \bigcap\limits_{i \in \omega} W_{x|_i} = \bigcap\limits_{i \in \omega} \bigcup\limits_{s \in 2^i} W_{s}
        \]
        is perfect. Clearly \( P \subseteq \bigcap\limits_{i \in \omega} U_i = A \). \(\)
    \end{Prf}

    It is worth noting that the multidimensional version of this proposition utilises basically the same proof idea, but it needs to be adjusted to the product topology.

    \begin{Prp}
        Let \( A \subseteq C^\omega \) be a dense \(G_\delta\) set. There exist non-empty perfect sets \( P_k \subseteq C \) for \( k \in \omega \) such that \( \prod\limits_{k \in \omega} P_k \subseteq A \).
    \end{Prp}
    \begin{Prf}
        Let \( A = \bigcap\limits_{i \in \omega} U_i  \), where the sets \( U_i \) are open and dense. Let
        \[
            V_0 = \prod\limits_{k \in \omega} V_{0,k} \subseteq U_0
        \]
        be a basic clopen set. Note that \( V_{0,k} \) is a proper subset of \( C \) for only a finite amount of \(k \in \omega\). We can represent \(V_{0,0}\) as a sum of two disjoint basic clopen sets \( W_{(0), 0}, W_{(1),0} \subseteq C \). Now suppose that we have a clopen set 
        \[
            V_i = \prod\limits_{k \in \omega} V_{i,k} \subseteq U_0 \cap \ldots \cap U_i
        \]
        such that for \( k \le i \) the sets \( V_{i,k} \) are represented as a disjoint union
        \[
            V_{i,k} = \bigcup\limits_{s \in 2^{i+1}} W_{s,k}
        \]
        of basic clopen subsets of \( C \), ie
        \[
            V_i = \bigcup\limits_{s \in (2^{i+1})^i} W_{s(0), 0} \times \ldots W_{s(i), i} \times \prod\limits_{k > i} V_{i,k}.
        \]
        Note that the sum above is finite. By taking intersections of \( U_{i+1} \) with all those sets one after another we can find the basic clopen sets \( W^*_{s,k} \subseteq W_{s,k} \) and for \( k > i \) the basic clopen sets \( V_{i+1,k} \subseteq V_{i,k} \) (of which only finite amount are proper subsets) such that the sum
        \[
            \bigcup\limits_{s \in (2^i)^i} W^*_{s(0), 0} \times \ldots W^*_{s(i), i} \times \prod\limits_{k > i} V_{i+1,k} \subseteq \subseteq U_0 \cap \ldots \cap U_{i+1}.
        \]
        We can represent each set \(W^*_{s,k}\) as a disjoint union of basic clopen sets \( W_{s^\smallfrown0,k}, W_{s^\smallfrown1} \subseteq C  \) and the set \( V_{i+1,i+1} \) as the sum of \( 2^{i+2} \) disjoint basic clopen sets \( W_{s, i+1} \) for \( s \in 2^{i+2} \). For \( k \le i \) we can take
        \[
            V_{i+1,k} = \bigcup\limits_{s \in 2^{i+2}} W_{s,k}.
        \]
        We consider the fusion sequence on each coordinate separately. From the properties of fusion we obtain that all the sets
        \[
            P_k = \bigcup\limits_{x \in 2^\omega} \bigcap\limits_{i \ge k} W_{x|_i,k} = \bigcap\limits_{i > k} \bigcup\limits_{s \in 2^i} W_{s,k}
        \]
        are perfect. Clearly we have \( \prod\limits_{k \in \omega} P_k \subseteq \bigcap\limits_{i \in \omega} U_i = A \).
    \end{Prf}

    We will now proceed with proving that we can substitute in the original problem functions from \( \R^\omega \) to \( [0;1] \) by functions from \( C^\omega \) to \( C \).
    
    \begin{Lem}
        Let \( P_k \subseteq C \) be perfect sets and let \( f \colon \prod\limits_{k \in \omega} P_k \to [0;1] \) be a continuous function. If the image of the function \(f\) is perfect then there exist perfect sets \( Q_k \subseteq P_k \), each homeomorphic to the Cantor set, such that \( f(\prod\limits_{k \in \omega}P_k) \) is either not perfect or is homeomorphic to the Cantor set.
    \end{Lem}
    \begin{Prf}
        Let \( a_{0,0} \) and \( a_{0,1} \) be the smallest and largest elements of the image of \( f \) respectively. We will define a sequence of refining partitions of \( [a_{0,0}; a_{0,1}] \) into interval in a following way. As the image of \( f \) is perfect clearly \( a_{0,0} < a_{0,1} \) and with \( a_{n,i} \) defined for \( i \le 2^n \) there have to exist \( a_{n+1,i} \in [0;1] \) for \( i \le 2^{n+1} \) such that
        \begin{itemize}
            \item \( a_{n+1, 2i} = a_{n, i} \) for \( i \le 2^n \)
            \item \( a_{n+1,i} < a_{n+1,i+1} \) for \( i < 2^{n+1} \)
            \item \( f^{-1}(a_{n+1,i}) \) is nowhere dense for \( 0 < i < 2^{n+1} \)
            \item \( f^{-1}([a_{n+1,0}; a_{n+1,1})), f^{-1}((a_{n+1,1}; a_{n+1,2})), \ldots ,\)
            
             \(f^{-1}((a_{n+1,2^{n+1}-2}; a_{n+1,2^{n+1} -1})) , f^{-1}((a_{n+1,2^{n+1}-1}; a_{n+1,2^{n+1}}]) \) are not nowhere dense
            \item the diameter of partitions defined in such a manner tends to \( 0 \) with increasing \(n\)
        \end{itemize}
        Let
        \[
            A_n = \{ a_{n,i} \colon i \le 2^n \}.
        \]
        Then we have \( A_{n} \subseteq A_{n+1} \). The set \( A = \bigcup_{n \in \omega} A_n \) is countable and its preimage \( f^{-1}(A) \) is a meager \( F_\sigma \) set, ie \( \prod\limits_{k \in \omega} P_k \setminus f^{-1}(A) \) is a dense \( G_\delta \) set a thus contains a product \( \prod\limits_{k \in \omega} Q_k \) where all the sets \( Q_k \) are perfect and non-empty. The image \( f(\prod\limits_{k \in \omega} Q_k) \) is clearly zero-dimensional and thus if it is perfect it has to be homeomorphic to the Cantor set.
    \end{Prf}

    \begin{Thm}
        Let the following property hold.
        \begin{itemize}
            \item For any continuous functions \( f_n \colon C^\omega \to C \) there exist a set \( N = \{ n_i \colon i \in \omega \} \in [\omega]^\omega \) and non-empty perfect sets \( P_k \subseteq C \) for \( k \in \omega \) such that the subsequence \( (f_{n_i})_{i \in \omega} \) is pointwise convergent on the product \( \prod\limits_{k \in \omega}  P_k \).
        \end{itemize}
        Then the answer to the Harrington problem is positive.
    \end{Thm}
    \begin{Prf}
        Let \( P_k \subseteq \R \) be perfect for \( k \in \omega \). Observe that if for an infinite amount of \( n \in \omega \) the set \( f_n(\prod\limits_{k \in \omega} P_k) \) is not perfect then we obtain the result right away as we can pick an isolated point in each of such sets and then find a convergent subsequence. Thus we can at any further point in the proof assume that all images \( f_n(\prod\limits_{i \in \omega} P_i) \) are perfect.
        
        We will begin with restricting all the functions to the set \( C^\omega \). By the lemma above there are perfect sets \( P_{0,k} \subseteq C \) for \( k \in \omega \) such that \( C_0 = f_0(\prod\limits_{k \in \omega} P_{0,k}) \) is homeomorphic to the Cantor set. Inductively we can define the sets \( P_{n+1,k} \subseteq P_{n,k} \) such that \( C_{n+1} = f_{n+1}(\prod\limits_{k \in \omega} P_{n+1,k}) \) is homeomorphic to the Cantor set. Moreover they can be chosen in such a way that
        \[
            P_{n,k} = \bigcup\limits_{s \in 2^{n+1}} Q_{s,k}
        \]
        for \( k \le n \), where all the sets \( Q_{s,k} \) are pairwise disjoint and they are intersections of \( P_{n,k} \) with basic clopen subsets in \( C \)
        
        Any finite sum of the sets \( C_i \) is perfect and zero-dimensional and thus homeomorphic to the Cantor set. Any infinite sum will be dense in itself, but it may not be closed. Let us pick such \( N \in [\omega]^{<\omega} \) that for any \( a,b \in [0;1] \) with \( a < b \) the set \( \bigcup\limits_{n \in N} C_n \cap (a;b) \) is not dense in \( (a;b) \). Then the closure \( D = \Cl(\bigcup\limits_{n \in N} C_n) \) is in fact homeomorphic to the Cantor set. From the properties of fusion we obtain that each set
        \[
            Q_k = \bigcup\limits_{x \in 2^\omega} \bigcap\limits_{i \ge k} W_{x|_i,k} = \bigcap\limits_{i > k} \bigcup\limits_{s \in 2^i} P_{s,k}
        \]
        contains a perfect subset \( P_k \). It follows that
        \[
            f_n|_{\prod\limits_{i \in \omega} P_i} \colon \prod\limits_{i \in \omega} P_i \to D
        \]
        and our assumption gives the desired result.
    \end{Prf}

    \section{Baire property and measurability}

    As one of our variants we consider Baire functions instead of continuous ones.
    
    \begin{Dfn}
        We will say that a subset \( A \) of a topological space \( X \) has Baire property if it can be represented as a symmetric difference \( U \triangle m \) of an open set \( U \) and a meager set \( m \).
        
        We will say that a function \( f \colon X \to Y \) - for topological spaces \( X \) and \( Y \) is Baire if preimage \( f^{-1}(U) \) of any open subset \( U \) of \( Y \) has Baire property. 
    \end{Dfn}

    A natural question arises if Baire functions are always continuous apart from some meager set. In \cite{KurBaire}, \cite{RFKun} and \cite{RFKul} one can find following partial answers to that question, which will be of importance to our considerations.
    \begin{Thm}
        Let \( X \) and \( Y \) be metric spaces. The following statements are equivalent:
        \begin{itemize}
            \item For every Baire function from \(X\) to \(Y\) there exists a meager set \( m \subseteq X \) such that \( f|_{X \setminus m} \) is continuous.
            \item There does not exist a partition (called \(K\)-partition) \( \mathcal{F} \) of \( X \) into meager subsets such that for any \( \mathcal{F}' \subseteq \mathcal{F} \) the sum \( \bigcup \mathcal{F}' \) has Baire property.
        \end{itemize}
    \end{Thm}
    \begin{Thm}
        Let \( f \colon X \to Y \) be Baire and \( Y \) be a separable metrizable space. Then there exists a meager set \( m \subseteq X \) such that \( f|_{X \setminus m} \) is continuous.
    \end{Thm}
    \begin{Thm}
        Let \( f \colon X \to Y \) be Baire and \( X \) be a completely metrizable space of weight at most \( \con \) and \( Y \) be a metrizable space. Then there exists a meager set \( m \subseteq X \) such that \( f|_{X \setminus m} \) is continuous.
    \end{Thm}

    As another one of our variants considers measurable functions the following variant of the well know Luzin's theorem, which can be found in \cite{Luz}, will be vital in our reasonings.
    
    \begin{Thm}
        (Luzin) Let \( E \subseteq \R \) be Lebesgue measurable and \( f \colon E \to R \). The function \( f \) is measurable iff for any \( \varepsilon > 0 \) there exists a closed set \( F \) such that \( f|_F \) is continuous and \( |E \setminus F| < \varepsilon \).
    \end{Thm}
    
    \section{Main result}
    
    We will apply fusion lemma to our problem.
    
     \begin{Thm}
        Let \( f_n \colon C \to 2 \) be continuous functions. Then there exist \( N = \{ n_i \colon i \in \omega \} \in [\omega]^\omega \) and a non-empty perfect set \( P \subseteq C \) such that the subsequence \( (f_{n_i})_{i \in \omega} \) is pointwise convergent on \(P\).
    \end{Thm}
    \begin{Prf}
        Let \( U_0^n = f_n^{-1}(0) \) and \( U_1^n = f_n^{-1}(1) \). Observe that those are disjoint clopen sets and their sum is whole \(C\). If infinitely many sets \( U_j^n \) are empty then the result follows in a straightforward way. Let \( n_0 \in \omega \) be such that both \( U_0^{n_0} \) and \( U_1^{n_0} \) are non-empty. There exist \( t_{(0)}, t_{(1)} \in 2^{<\omega} \) such that \[ U_{(0)} = C_{t_{(0)}} \subseteq U_0^0 \] and \[ U_{(1)} = C_{t_{(1)}} \subseteq U_1^0. \]
        
        Assume that all \( U_s \) and \( n_i \) are defined for \( s \in 2^{m+1} \) and \( i \le m \). If there exists \( s \in 2^{m+1} \) such that \( U_0^n \cap U_s = \emptyset \) or \( U_1^n \cap U_s = \emptyset \) for infinitely many \( n > n_m \) then once again the result follows in a straightforward way. Otherwise let \( n_{m+1} \in \omega \) be such that \( U_j^{n_{m+1}} \cap U_s \not= \emptyset \) for all \( s \in 2^{m+1} \) and \( j \in 2 \). It follows that there exist \( t_{s^\smallfrown 0}, t_{s^\smallfrown 1} \supset t_s \) such that \[ U_{s^\smallfrown 0} = C_{t_{s^\smallfrown 0}} \subseteq U_0^{n_{m+1}} \cap U_s \] and \[ U_{s^\smallfrown 1} = C_{t_{s^\smallfrown 1}} \subseteq U_1^{n_{m+1}} \cap U_s. \] Now consider a set \[ S = \{ x \in 2^{\omega} \colon x(2m) = 0 \mbox{ for } m \in \omega \}\] it is clearly uncountable. By the properties of fusion we obtain that the set
        \[
            Q = \bigcup\limits_{x \in S} \bigcap\limits_{i \in \omega} U_{x|_i} = \bigcap\limits_{m \in \omega} \bigcup\limits_{s \in 2^{2m+1}} U_{s^\smallfrown 0}
        \]
        is compact as well as uncountable and thus contains a perfect subset \( P \). It is easy to see that the sequence \( (f_{n_{2i}})_{i \in \omega} \) is convergent to \( 0 \) on \(P\).
    \end{Prf}

    \begin{Thm}
        Let \( f_n \colon C^\omega \to 2 \) be continuous functions. Then there exist \( N = \{ n_i \colon i \in \omega \} \in [\omega]^\omega \) and non-empty perfect sets \( P_k \subseteq C \) for \( k \in \omega \) such that the subsequence \( (f_{n_i})_{i \in \omega} \) is pointwise convergent on the product \( \prod\limits_{k \in \omega} P_k \). What is more we can assume that \( P_k = C \) for \( i > 0 \).
    \end{Thm}
    \begin{Prf}
        For any fixed \( x \in C^\omega \) let us define functions \( g_{x,n} \colon C \to 2 \) in a following way
        \[
            g_{x,n}(y) = f_n(y, x_0, x_1, \ldots).
        \]
        Let \( A = \{ a_m \colon m \in \omega \} \) be a countable, dense subset of \( C^\omega \). From the theorem above there exists a non-empty perfect set \( Q_\emptyset \subseteq C \) and \( N_0 \in [\omega]^\omega \) such that the subsequence \( (g_{a_0,n})_{n \in N_0} \) is pointwise convergent on \( Q_\emptyset \). There exist disjoint clopen sets \( U_{(0)}, U_{(1)} \subseteq C \) such that
        \[
            P_{(0)} = U_{(0)} \cap Q_\emptyset \not= \emptyset
        \]
        and
        \[
            P_{(1)} = U_{(1)} \cap Q_\emptyset \not= \emptyset.
        \]
        The sets \(P_{(0)}\) and \(P_{(1)}\) are perfect and thus homeomorphic to \(C\) itself.
        
        Assume that all \( P_s \) and \( N_i \) are defined for \( s \in 2^{m+1} \) and \( i \le m \). As the set \( 2^{m+1} \) is finite from the theorem above we get that there exist the perfect sets \( Q_s \subseteq P_s \) for \( s \in 2^{m+1} \) and \( N_{m+1} \in [N_m]^\omega \) such that the sequence \( (g_{a_{m+1},n})_{n \in N_{m+1}} \) is pointwise convergent on all the sets \( Q_s \). For each such set there exist disjoint clopen sets \( U_{s^\smallfrown 0}, U_{s^\smallfrown 1} \subseteq C \) such that
        \[
            P_{s^\smallfrown 0} = U_{s^\smallfrown 0} \cap Q_s \not= \emptyset
        \]
        and
        \[
            P_{s^\smallfrown 1} = U_{s^\smallfrown 1} \cap Q_s \not= \emptyset.
        \]
        Clearly all the sets \( P_{s^\smallfrown j} \) are homeomorphic to \(C\).
        
        Let us define the set \( N = \{ n_m \colon m \in \omega \} \) in a following way.
        \[
            n_0 = \min (N_0)
        \]
        and
        \[
            n_{m+1} = \min (N_{m+1} \setminus \{ n_0, \ldots, n_m \} ).
        \]
        Clearly \( (g_{a_m,n})_{n \in N} \) is convergent on all \( P_s \) for \( m \in \omega \) and \( s \in 2^{<\omega} \). By the properties of fusion we obtain that the set
        \[
            Q = \bigcup\limits_{x \in 2^\omega} \bigcap\limits_{i \in \omega} P_{x|_i} = \bigcap\limits_{m \in \omega} \bigcup\limits_{s \in 2^m} P_{s}
        \]
        is compact as well as uncountable and thus contains a perfect subset \( P \). We get that  \( (g_{a_m,n})_{n \in N} \) is convergent on \(P\) for any \( m \in \omega \). Thus it follows from density of \(A\) and continuity of the functions \(f_n\) we obtain that the sequence \( (f_n)_{n \in N} \) is pointwise convergent on the product \( P \times \prod\limits_{k > 0} C \).
    \end{Prf}

    \begin{Cor}
        Let \( f_n \colon C^\omega \to C \) be continuous functions. Then there exists \( N = \{ n_i \colon i \in \omega \} \in [\omega]^\omega \) and non-empty perfect sets \( P_k \subseteq C \) for \( k \in \omega \) such that the subsequence \( (f_{n_i})_{i \in \omega} \) is pointwise convergent on the product \( \prod\limits_{k \in \omega} P_k \). What is more we can assume that \( P_k = C \) for \( i > 0 \).
    \end{Cor}
    \begin{Prf}
        For \( x = (x_0, x_1, \ldots) \in 2^\omega = C \) define
        \[
            \pi_l(x) = x_l \mbox{ for } l \in \omega.
        \]
        The projections \(pi_n\) are clearly continuous. Let \( g_{l,n} = \pi_l \circ f_n \). We can apply the theorem above to the functions \( g_{0,n} \) and get the non-empty perfect set \( Q_\emptyset \subseteq C \) and \( N_0 \in [\omega]^\omega \)  such that \( (g_{0,n})_{n \in N} \) is pointwise convergent on \( Q_\emptyset \times \prod\limits_{k > 0} C \). There exist disjoint clopen sets \( U_{(0)}, U_{(1)} \)  such that
        \[
            P_{(0)} = U_{(0)} \cap Q_\emptyset \not= \emptyset
        \]
        and
        \[
            P_{(1)} = U_{(1)} \cap Q_\emptyset \not= \emptyset.
        \]
        The sets \(P_{(0)}\) and \(P_{(1)}\) are perfect and thus homeomorphic to \(C\) itself.
        
        Assume that all \( P_s \) and \( N_i \) are defined for \( s \in 2^{l+1} \) and \( i \le l \). As the set \( 2^{l+1} \) is finite from the theorem above we get that there exist the perfect sets \( Q_s \subseteq P_s \) for \( s \in 2^{l+1} \) and \( N_{l+1} \in [N_l]^\omega \) such that the sequence \( (g_{l+1,n})_{n \in N_{l+1}} \) is pointwise convergent on the product \( (\bigcup\limits_{s \in 2^{l+1}} Q_s) \times \prod\limits_{i > 0} C \). For each such set there exist disjoint clopen sets \( U_{s^\smallfrown 0}, U_{s^\smallfrown 1} \subseteq C \) such that
        \[
            P_{s^\smallfrown 0} = U_{s^\smallfrown 0} \cap Q_s \not= \emptyset
        \]
        and
        \[
            P_{s^\smallfrown 1} = U_{s^\smallfrown 1} \cap Q_s \not= \emptyset.
        \]
        Clearly all the sets \( P_{s^\smallfrown j} \) are homeomorphic to \(C\).
        
        Let us define the set \( N = \{ n_m \colon m \in \omega \} \) in a following way.
        \[
            n_0 = \min (N_0)
        \]
        and
        \[
            n_{m+1} = \min (N_{m+1} \setminus \{ n_0, \ldots, n_m \} ).
        \]
        By the properties of fusion we obtain that the set
        \[
            Q = \bigcup\limits_{x \in 2^\omega} \bigcap\limits_{i \in \omega} P_{x|_i} = \bigcap\limits_{l \in \omega} \bigcup\limits_{s \in 2^l} P_{s}
        \]
        is compact as well as uncountable and thus contains a perfect subset \( P \). We get that  \( (g_{l,n})_{n \in N} \) is convergent on \( P \times \prod\limits_{i>0} C \) for any \( l \in \omega \).
    \end{Prf}

    Applying our earlier codomain reduction argument we obtain.
    
    \begin{Cor}
        Let \( X_k \) be metric spaces and \( Q_k \subseteq X_k \) be perfect for \( k \in \omega \). For any continuous functions \( f_n \colon X_k \to [0;1] \) there exist \( N = \{ n_i \colon i \in \omega \} \in [\omega]^\omega \) and non-empty perfect sets \( P_k \subseteq Q_k \) for \( k \in \omega \) such that the subsequence \( (f_{n_i})_{i \in \omega} \) is pointwise convergent on the product \( \prod\limits_{k \in \omega} P_k \).
    \end{Cor}

    Amongst other this result provides a positive answer to the original Harrington problem. The generalization to the Baire functions follows in a straightforward way.

    \begin{Thm}
        Let \( f_n \colon C^\omega \to [0;1]\) be Baire functions. Then there exists \( N = \{ n_i \colon i \in \omega \} \in [\omega]^\omega \) and non-empty perfect sets \( P_k \subseteq C \) for \( k \in \omega \) such that the subsequence \( (f_{n_i})_{i \in \omega} \) is pointwise convergent on the product \( \prod\limits_{k \in \omega} P_k \).
    \end{Thm}
    \begin{Prf}
        As weight of \( C^\omega \) is equal to \( \omega \) it follows that each \( f_n \) is continuous apart from a meager set, ie it is continuous on the intersection
        \[
            G_n = \bigcap\limits_{i \in \omega} U_{n,i}
        \]
        of open and dense subsets of \( C^\omega \). As \( C^\omega \) is a Baire space the set
        \[
            G = \bigcap\limits_{n \in \omega} G_n = \bigcap\limits_{(n,i) \in \omega^2} U_{n,i} = \bigcap\limits_{j \in \omega} U_j
        \]
        is a dense \( G_\delta \) set and all of the functions \( f_n \) are continuous on \(G\). We obtain \( P = \prod\limits_{k \in \omega} P_k \subseteq G \) such that all \( P_k \) are homeomorphic to the Cantor set. As all the functions \( f_n \) are continuous on \( G \) they are also continuous on \( P \) and the result follows directly from the theorems above.   
    \end{Prf}

    \section{Further developments}
    
    In the next part we will generalize our results to the wider variety of functions as well as topological spaces, including
    \begin{enumerate}
        \item measurable functions
        \item functions with \((s)\)-property which are modelled after Sack's forcing \cite{RFAni}
        \item functions with an analog of \( (s) \)-property for Silver's forcing
        \item completely Ramsey function on the space \( [\omega]^{<\omega} \) with the Ellentuck topology \cite{El} (which could be thought of as a topological representation of Mathias forcing)
    \end{enumerate}
    Applying fusion technique to those cases will remain our main focus. Some of those variants might require using fusion technique for different forcing notions and the generalizations of fusion such as Axiom A. \cite{BaugMat}

\end{document}